\documentclass[12pt]{article}
    \usepackage{amsmath, amssymb, amsthm, amsxtra,    }

\def\be{\begin{equation}}
\def\ee{\end{equation}}

\def\vp{\varphi}
\def\arrowk{^\to{\kern -6pt\topsmash k}}
\def\arrowK{^{^\to}{\kern -9pt\topsmash K}}
\def\arrowt{^\to{\kern -6pt\topsmash t}}
\def\arrowr{^\to{\kern-6pt\topsmash r}}
\def\arrowvp{^\to{\kern -8pt\topsmash\vp}}
\def\tk{\tilde{\kern 1 pt\topsmash k}}
\def\barm{\bar{\kern-.2pt\bar m}}
\def\barN{\bar{\kern-1pt\bar N}}
\def\barA{\, \bar{\kern-3pt \bar A}}

\def\be{\begin{equation}}
\def\ee{\end{equation}}
\numberwithin{equation}{section}
\begin{document}
    \theoremstyle{plain}
    \newtheorem{theorem}{Theorem}
    \newtheorem{lemma}{Lemma}
    \newtheorem{proposition}{Proposition}
    \newtheorem{corollary}[theorem]{Corollary}
    \theoremstyle{definition}
    \newtheorem*{definition}{Definition}
    \theoremstyle{remark}
    \newtheorem*{remark}{Remark}

\title{ A Remark on Sieving in Biased Coin Convolutions }
\author{Mei-Chu Chang\footnote{Research partially
financed by the NSF Grants~DMS~1301608.}\\ \texttt{Department of Mathematics}\\
\texttt{University of California, Riverside}\\\texttt{\small
mcc@math.ucr.edu}}

\date{}
\maketitle
\begin{center}{\bf Abstract}\end{center}

{\it In this work, we establish a nontrivial level of distribution for densities on $\{1,\ldots, N\}$ obtained by a biased coin convolution. As a consequence of sieving theory, one then derives the expected lower bound for the weight of such densities on sets of pseudo-primes. }

\bigskip

\noindent {\bf \large Introduction.}

Over the recent years, there has been an increasing interest in sieving problems in combinatorial objects without a simple arithmetic structure. The typical example is that of finitely generated `thin subgroups' of linear groups such as $SL_2(\mathbb Z)$ or $SL_2(\mathbb Z + i\mathbb Z)$. These groups are combinatorially defined but are not arithmetic (they are of infinite index) and as such cannot be studied with classical automorphic techniques. Examples of natural appearances of this type of questions include the study of the curvatures in integral Apollonian circle packings, Pythagorean triples and issues around fundamental discriminates of quadratic number fields and low lying geodesics in the modular surface. (See \cite{BK}.)  The reader may also wish to consult the excellent Bourbaki exposition by E. Kowalski \cite{K} for a detailed account of many of these recent developments around `exotic sieving'.

In this paper we consider a slightly different problem but in a somewhat similar spirit. Let $N=2^m$ and identify $\{1,\ldots, N\}$ with the Boolean cube $\{0,1\}^m$ through binary expansion. Denote $\mu_{\rho}$ the probability measure on $\{0,1\}^m$ given by a standard biased coin convolution, i.e. on each factor we take an independent distribution assigning probability $\rho$ to $0$ and $1-\rho$ to $1$. Consider the resulting distribution on $\{1,\ldots, N\}$. For $\rho=\frac 12$, this is the uniform distribution while for $\rho \to 1$, these distributions become increasingly singular. Our aim is to study some of their arithmetical properties and in particular prove that there is a nontrivial level of distribution no matter how close $\rho$ is to $1$, $\rho <1$. Similar results may also be obtained for $g$-adic analogues, expanding integers in base $g$. 

\bigskip

\noindent{\bf\large Notations.}

\medskip

$e(\theta)=e^{2\pi i\theta}$, $e_q(\theta)=e(\frac {\theta}q).$

$ c, C=$ various constants.

$A\ll B$ and $A=O(B)$ are each equivalent to that $|A|\leq cB$ for some constant $c$. If the constant $c$ depends on a parameter $\rho$, we use $\ll_\rho$. Otherwise, $c$ is absolute.

\section{The statement.}

Consider the distribution $\mu$ on $[1, N]\cap\mathbb Z$, with $N=2^m$, induced by the random variable $\sum_j \xi_j 2^j$
with $(\xi_j), j\ge 0$, be an independent, identically distributed sequence of random variables taking values in $\{0, 1\}$, $\mathbb P[\xi_j=0]=\rho$, $\mathbb P[\xi_j=1]=1-\rho$, $\frac 12<\rho<1$. Thus, if $n=\sum_ja_j2^j$ with $a_j\in\{0,1\}$ the binary expansion, then  \be\label{1.1}\mu(n)=\rho^{m-\ell}(1-\rho)^\ell, \;\;\text{ where } \ell=\sum_ja_j\ee Note that for $\rho=\frac 12$ we obtain the normalized uniform measure on $[0, N]$.

The measure \eqref{1.1} has dimension $(1-\rho)\log\frac 1{1-\rho}$ and hence becomes more irregular for $\rho\to 1$. Our aim is to establish a level of distribution of $\mu$ in the sense of sieving theory. Thus, taking $q<N^{\alpha}$, $q$ square free and $\alpha$ appropriately small, (since $\mu$ is normalized) we may write
\be\label{1.2}\begin{aligned}\mu\big[n\le N:q|n\big]=&\frac 1q\sum_{\lambda=0}^{q-1}\sum_{n=1}^N e_q(\lambda n)\mu(n)\\
=&\frac 1q+R_q,\end{aligned}\ee
where $$R_q=\frac 1q\sum_{\lambda=1}^{q-1}\sum_{n=1}^N e_q(\lambda n)\mu(n).$$
We also assume $q$ odd. The number $\alpha$ is the {\it sieving exponent}.

Our aim is to obtain a bound of the form
\be\label{1.3}
\sideset{}{'}\sum_{q<N^{\alpha}}|R_q|= o(1)\ee
where $\sum'$ sums over $q$ square free and odd.

\medskip
\noindent{\bf Theorem 1.} {\it Let the notations be as above. Then $\mu$ has sieving exponent $\alpha(\rho)>0$. In fact, $\alpha(\rho)=O(1-\rho)$ for $\rho\to 1$.}

\medskip

Sieving pseudo primes is the goal of sieving theory. From standard combinatorial sieve (which also applies to measures instead of sets.) (See e.g. \cite{BGS}, \cite{BK}, \cite{FI}, \cite{IK}) we have the following result about {\it r-pseudo-primes} (products of at most r primes).

\medskip

\noindent{\bf Corollary 2.} {\it
\be\label{1.4}\mu(\mathcal P_r\cap [0, N])\sim \frac 1{\log N}\ee
with $\mathcal P_r=\{ r\text{-pseudo-primes}\}, r=r(\rho)$.}

\bigskip
\section{First estimates.}

\medskip

Let
\be\label{2.1}
\begin{aligned}
R_q=&\frac 1q\sum_{\lambda=1}^{q-1}\sum_{n=1}^N e_q(\lambda n)\mu(n)\\
=&\frac 1q\sum_{\lambda=1}^{q-1}\prod_{j<m}\bigg(\rho+(1-\rho)\;e\bigg(\frac{\lambda 2^j}q\bigg)\bigg).\end{aligned}\ee
Note that
\be\label{2.2}
|\rho+(1-\rho)e(\theta)|^2=1-4\rho(1-\rho)\sin^2 \pi \theta .\ee
Let us consider first the case of small $q$.

For $\lambda\not=0\mod q$, \eqref{2.2} implies
$$\bigg|\rho+(1-\rho)\;e\bigg(\frac{\lambda 2^j}q\bigg)\bigg|\le 1-\frac c{q^2}$$ for $c>0$
so that \eqref{2.1}$<\big(1-O(\frac 1{q^2})\big)^m< e^{-C\frac m{q^2}}<N^{-c/q^2}$.

One can do better by the following observation.

Let $I\subset\{1,\ldots,m\}$ be an arbitrary interval of size $\sim \log q$. Then for $\lambda\not=0\mod q$,
\be\label{2.3}
\max\bigg\{\sin^2\frac {\lambda 2^j}q \pi : j\in I\bigg\}>c\ee
with $c>0$ some constant independent of $q$. Therefore, we also have
\be\label{2.4}
(2.1)< \big( 1-c(\rho)\big)^{\frac m{\log q}} < N^{-\frac {c(\rho)}{\log q}} < e^{-\sqrt {\log N}}\ee
if $\log q< O\big(\sqrt{\log N}\big)$.

\bigskip
\section{Further estimates.}
 We want to estimate
 \be\label{3.1}
 \sum_{q\sim Q}|R_q|
 \ee
 with $Q<N^{\alpha}$ and $\log Q \gtrsim \sqrt{\log N}$. It will suffice to show that (3.1)$<Q^{-c}$ for some $c>0$.

 We may assume $\alpha=\frac 1t$ for some large $t\in \mathbb Z$ \big(given in \eqref{3.7}\big). Choose $h\in \mathbb Z$ such that
 \be\label{3.2}
 2^h\sim Q^2.\ee
 Hence
 $$h<\frac 2tm<m.$$
 Estimate \eqref{3.1} using H\"{o}lder inequality
 \be\label{3.3}
 \begin{aligned}
 &\sideset{}{'}\sum_{q\sim Q}|R_q|\\
 \le \;\; &\sideset{}{'}\sum_{q\sim Q}\;\frac 1Q\;\sum_{\lambda=1}^{q-1}\;\prod_{\tau=1}^{t/2}\;\;\prod_{j=(\tau-1)h}^{\tau h}
 \bigg|\rho+(1-\rho)\;e\bigg(\frac{\lambda 2^j}q\bigg)\bigg|\\
 \le \;\;\;\;&\sideset{}{'}\sum_{q\sim Q}\;\bigg[\prod_{\tau=1}^{t/2}\frac 1Q \;\sum_{\lambda=1}^{q-1} \;\prod_{j=(\tau-1)h}^{\tau h}
 \bigg|\rho+(1-\rho)\;e\bigg(\frac{\lambda 2^j}q\bigg)\bigg|^{t/2}\bigg]^{2/t}\\
 =\;\;\;\;&\sideset{}{'}\sum_{q\sim Q}\frac 1Q\;\sum_{\lambda=1}^{q-1} \;\prod_{j=0}^{h-1}
 \bigg|\rho+(1-\rho)\;e\bigg(\frac{\lambda 2^j}q\bigg)\bigg|^{t/2}.
 \end{aligned}
 \ee
 For the last equality, we note that for each $\tau$
 $$\begin{aligned}&\{\lambda 2^j\text{ mod } p: (\tau-1)h\le j<\tau h \}\\
 =&\{\lambda 2^j\text{ mod } p: 0\le j< h \}.\end{aligned}$$

 \medskip
To finish the estimate, we need the following two lemmas.

 \medskip

 \noindent{\bf Lemma 3.} {\it For all $\theta$, $0<\delta<1$ and
 \be\label{3.4}
 \ell >\frac {\log{\frac 1{\delta}}}{\rho(1-\rho)},
 \ee
 we have
 \be\label{3.5}|\rho+(1-\rho)e(\theta)|^{2\ell}\le 1-(1-\delta)\sin^2 \pi\theta.\ee
 }

 \medskip

 \noindent{\it Proof.} Let
 $$\gamma= 4\rho(1-\rho)\sin^2 \pi \theta .$$
 By \eqref{2.2},
 $$|\rho+(1-\rho)e(\theta)|^{2\ell}=1-\gamma.$$
 We consider the following two cases.

 \smallskip

\noindent (i). $\gamma>\frac 1{\ell}\log\frac 1{\delta}$.

 \smallskip

 Then
 $$(1-\gamma)^{\ell} \le e^{-\ell\gamma}<\delta<1-(1-\delta)\sin^2 \pi\theta.$$

 \smallskip

 \noindent(ii). $\gamma\le\frac 1{\ell}\log\frac 1{\delta}$.

 \smallskip

 Let
 $$\ell_1 =\frac{\ell}{2\log\frac 1{\delta}}<\ell$$
 and estimate
 $$\begin{aligned}(1-\gamma)^{\ell}< (1-\gamma)^{\ell_1}<e^{-\ell_1\gamma}&< 1-\frac 12\ell_1\gamma\\
 &= 1- \frac{\ell\rho(1-\rho)}{\log\frac 1{\delta}}\sin^2 \pi\theta\\
 &< 1-\sin^2 \pi\theta\\
 &< 1-(1-\delta)\sin^2 \pi\theta.
 \end{aligned}$$
 (Note that the third inequality is because $\ell_1\gamma < \frac 12$.) $\;\; \square$

 \medskip

 \noindent{\bf Lemma 4.} {\it Let $\gamma< 1/10$ be positive. Then for all $\theta$ and $0<\delta<1$, we have
 \be\label{3.6}
 1-(1-\delta)\sin^2 \theta \le 1+\gamma-(1-\delta)\sin^2 (\theta+\gamma).
 \ee
 }

 \medskip

 \noindent{\it Proof.} Using the identity
 $$\sin^2A-\sin^2B=\sin(A+B)\sin(A-B)$$
 on the difference of both sides of \eqref{3.6}, we obtain
 $$(1-\delta)\big(\sin(2\theta+\gamma)\;\sin \gamma\big),$$
 which is bounded by $\gamma.\quad\square$

 \medskip

 Let \be\label{3.7}t>\frac {4\log{\frac 1{\delta}}}{\rho(1-\rho)}.\ee With $\theta=\lambda 2^j/q$, Lemma 3 implies that \eqref{3.3} is bounded by
 \be\label{3.8}
  \frac 1Q\sideset{}{'}\sum_{q\sim Q}\;\sum_{\lambda=1}^{q-1} \;\prod_{j=0}^{h-1}
 \bigg(1-(1-\delta)\sin^2\big(\frac{\pi \lambda 2^j}q\big) \bigg).
  \ee
    Given $Q$, let
  $$S=\bigg\{\frac {\lambda}q: 0\le \lambda <q, \;\;q\sim Q\bigg\}\subset [0,1].$$
  We note that $|S|\sim Q^2$ and $S$ is $Q^{-2}\sim 2^{-h}\;$ separated.

         In Lemma 4, taking $\gamma= \pi 2^j\beta'$ with $\beta'\in [0, \beta]$ for some $\beta=O(2^{-h})$ to be specified later, we bound \eqref{3.8}  by
  \be\label{3.9}\frac 1Q\sum_{\frac {\lambda}q\in S} \;\prod_{j=0}^{h-1}
 \bigg(1+\gamma-(1-\delta)\sin^2 \big(\pi  2^j(\frac{\lambda}q+\beta')\big) \bigg)\ee

  We will use integration to bound \eqref{3.9} by replacing $S$ by $S_{\beta}=S+[0, \beta]$. Averaging over $\beta'\in[0, \beta]$ gives
  \be\label{3.10}\begin{aligned}
  &\frac 1{\beta Q}\int_{S_{\beta}}\;\prod_{j=0}^{h-1}
 \big(1+\gamma-(1-\delta)\sin^2 (\pi 2^j x) \big)dx\\
 \lesssim &\frac 1{\beta Q}\int_0^1\;\prod_{j=0}^{h-1}
 \big(1+\gamma-(1-\delta)\sin^2 (\pi 2^j x) \big)dx
 \end{aligned}.\ee
  More precisely, we take
   \be\label{3.11}
   \beta=\frac{\;\delta}4\; Q^{-2},
    \ee
    (which implies $\gamma<\delta$) and bound \eqref{3.10} by
    \be\label{3.12}
    \begin{aligned}
    &\frac 4{\delta}\;Q\int_0^1 \;\prod_{j=0}^{h-1}
 \big(1+\delta-(1-\delta)\sin^2 (\pi 2^j x) \big)dx\\
 =\;&\frac 4{\delta}\;Q\bigg(1+\delta-\frac{1-\delta}2\bigg)^h\\
 =\; &\frac 4{\delta}\;Q\bigg(\frac{1+3\delta}2\bigg)^h\\
 <\;&Q^{-1/2},
 \end{aligned}
 \ee
 for $\delta$ small enough.

 Putting \eqref{3.3}, \eqref{3.8}-\eqref{3.10} and \eqref{3.12} together, we obtain the intended bound on \eqref{3.1}.

 \bigskip

 \section{Random polynomials with coefficients in $\{0, 1, -1\}$.}

The initial motivation for this work came from  \cite{PSZ}, where 
 one considers biased coin convolution densities for ternary expansions, with probabilities $\mathbb P[\xi=0]=\rho_0$, $\mathbb P[\xi=1]=\rho_1$, $\mathbb P[\xi=-1]=\rho_{-1}$ and $\rho_0 \ge \rho_1, \rho_{-1}$. The main problem focused in \cite{PSZ}  is to ensure that the set of integers $\{n<N: q^2|n \text{ for some }  q> Q\}$ carries small weight for $Q\to \infty$, which they manage to ensure if $q$ is not too large. The natural problem is whether such restriction is necessary. Clearly, this issue may be rephrased as the sieving problem for square free integers, but with unrestricted level of distribution. (The large values of $q$ are indeed the problematic ones.) While we are unable to provide a definite answer to their question and the main result of this note does not directly contribute, we will point out a simple probabilistic argument leading to the replacement of their condition.
Our argument uses virtually no arithmetic structure. 

Let $(\xi_j), j\ge 0$, be an independent, identically distributed sequence of random variables taking values in $\{-1, 0, 1\}$. Let $m\ge 1$ and define the random polynomial $P$ by
$$P(z):= \sum_{j=0}^m\xi_jz^j$$
In \cite{PSZ}, the authors assumed that
\be\label{1}
\max_{x\in\{-1, 0, 1\}}\mathbb P(\xi_0=x)< \frac 1{\sqrt 3}=0.5773\ldots
\ee
and proved that $\mathbb P(P \text{ has a double root }) = \mathbb P(P \text{ has } -1, 0 \text{ or } 1 \text{ as a double root })$ up to a $o(m^{-2})$ factor, and $\lim_{m\to \infty}\mathbb P(P \text{ has a double root }) =\mathbb P(\xi_0=0)^2$. One of the open problems they raised at the end of the paper asked whether it is necessary to have assumption \eqref{1}, which enters into the proof mainly through Claim 2.2 in their paper (which is crucial to their results). In this note, we will prove Claim 2.2 under a weaker assumption than assumption \eqref{1}. More precisely, we prove the following.

\bigskip

\noindent{\it
Assume
\be\label{1.1.1}
\max_{x\in\{-1, 0, 1\}}\mathbb P(\xi_0=x)< 0.7615\ldots.
\ee
Then there exist constants $C, c>0$ such that for any $B>0$ we have
\be\label{2}
\mathbb P(P(3) \text{ is divisible by } k^2 \text{ for some } k\ge B) \le CB^{-c}.
\ee}

\smallskip

\noindent{\it Remark.} The bound in \eqref{1.1.1} is the solution to equation \eqref{7}.

\medskip

 \noindent{\bf Proof.} Fix $r$ such that
 \be\label{2}
 3^r\le B^2< 3^{r+1}.
 \ee

 \noindent
 {\it Claim.
\be\label{3}
\mathbb P(P(3) \text{ is divisible by } k^2 \text{ for some } k \in[B, 2B]) \le 2^{-cr}
\ee
for some constant $c>0$.}

\medskip

\noindent
 {\bf Proof of Claim.} We write
 $$P(3)=\sum_{j<r}\xi_j3^j+\sum_{j=r}^m\xi_j3^j.$$
 Fix $\xi_r, \ldots, \xi_n$, and let $\ell= \sum_{j=r}^m\xi_j3^j.$

 If $k^2$ divides $P(3)$, then
 $$\sum_{j<r}\xi_j3^j\equiv -\ell \text{ mod } k^2.$$
 Since $\big|\sum_{j<r}\xi_j3^j\big|< 3^r/2 \le k^2/2$, we may denote $$\ell(k):= \sum_{j<r}\xi_j3^j\in \bigg(\frac{-k^2}2, \frac{k^2}2\bigg)$$
 and let $$S=\big\{\ell(k): k\in [B, 2B]\big\}\subset \big(-2B^2, 2B^2\big).$$
 It follows that
 \be\label{4}
 \text{ the left-hand-side of } (4.5) \le \mathbb P\big(\sum_{j<r}\xi_j3^j \in S\big).\ee

Let $\sigma_{(k)}= \big(\sigma_{(k)}(j)\big)_{j=0,\ldots, r-1}\in\{-1,0,1\}^r$ be defined by
$$\sum_{j<r}\sigma_{(k)}(j)3^j=\ell(k)$$
and let
$$A=\{\sigma_{(k)}: k\in[B, 2B]\} \text{ with } |A|\sim \sqrt{3}^{\;r}.$$
Let $\delta_j$ be the indicator function of $j, j=-1, 0,1$, and denote
$$\rho_j:=\mathbb P(\xi_0=j)\;\text{ for } \; j=-1, 0, 1, \;\;\text{ and }\;\; \rho:=\max_j \rho_j.$$
Denote the product measure on $\{-1,0,1\}^r$ by
$$\nu:=\bigotimes_{j=0}^{r-1} (\rho_0\delta_0+\rho_1\delta_1+\rho_{-1}\delta_{-1}).$$ Therefore we have
 (reasoning given below the display)
\be\label{11}\begin{aligned}(4.6)\le &\sum_{\sigma\in A} \nu(\sigma)\\
\le &|A|^{ 1/p}\bigg(\sum_{\sigma\in A} \nu(\sigma)^q\bigg)^{1/q}, \quad\text{with } \;\; \frac 1p +\frac 1q=1\\
\lesssim & \sqrt 3\;^{r/p} \big( \rho_0^q + \rho^q_1 + \rho^q _{-1}\big)^{r/q}\\
\le & \sqrt 3\;^{r/p} \big( \rho^q + (1-\rho)^q\big)^{r/q}\\
< & 2^{-cr}  \quad\text{for some constant }\;\; c>0.
\end{aligned}\ee
The second inequality is by H\"{o}lder, and the third inequality follows from the following estimate.
$$\begin{aligned}
\sum_{\sigma\in A} \nu(\sigma)^q=&\sum_{\sigma\in A} \bigotimes_{j=0}^{r-1} \big(\rho_0\delta_0(\sigma(j))+\rho_1\delta_1(\sigma(j))+\rho_{-1}\delta_{-1}(\sigma(j))\big)^q\\
=&\sum_{\sigma\in A} \bigotimes_{j=0}^{r-1} \big(\rho_0^q\delta_0(\sigma(j))+\rho_1^q\delta_1(\sigma(j))+\rho_{-1}^q\delta_{-1}(\sigma(j))\big)\\
\le & \sum_{a+b+c=r} \binom{r}{a}\binom{r-a}{b}\rho_0^{aq}\rho_1^{bq}\rho_{-1}^{cq}=(\rho_0^q+\rho_1^q+\rho_{-1}^q)^r.
\end{aligned}$$
To finish the proof of the claim, we want to show \eqref{11} $< 2^{-cr}$ for some constant $c>0$, i.e.
 $$\sqrt 3\;^{1/p} \big( \rho^q + (1-\rho)^q\big)^{1/q}<1,$$
and we want to solve
\be\label{5}t^q + (1-t)^q= \bigg(\frac 1{\sqrt 3}\bigg)^{\frac 1{p-1}},  \quad\text{with } \;\; \frac 1p +\frac 1q=1.\ee
Let $u=\frac 1{p-1}$ and rewrite \eqref{5} as
\be\label{6}\big( t^{1+u}+(1-t)^{1+u}\big)^{1/u}=\frac 1{\sqrt 3}\ee
Let $p$ go to infinity (hence $u$ goes to $0$). Then
$$ \begin{aligned}&t^{1+u}+(1-t)^{1+u}\\
= \;& t(1+u\log t +O(u^2))+  (1-t)(1+u\log (1- t) +O(u^2))\\
= \;& 1+\big(t\log t +(1-t)\log(1-t)\big)u +O(u^2).
\end{aligned}
$$
Hence \eqref{6} becomes
$$\bigg( 1+\big(t\log t +(1-t)\log(1-t)\big)u +O(u^2) \bigg)^{1/u}=\frac 1{\sqrt 3}.$$
In the limit for $u\to 0$, we obtain
$$e^{t\log t+(1-t)\log(1-t)}=\frac 1{\sqrt 3}.$$
Solving
\be\label{7}t^t(1-t)^{1-t}=\frac 1{\sqrt 3},\ee
we obtain $t=0.7615332817632392\cdots. \quad\square$

\medskip

It is possible to exploit somewhat better arithmetical features of the distribution under considerations but gains turn out to be minimal ($0.7654$ from $0.7615$), therefore, will not be elaborated here.

\bigskip

\bigskip

\bigskip

\bigskip
\noindent{\it Acknowledgement.}  The author would like to thank Gwoho Liu for computer assistance.

\end{document}